\documentclass[12pt]{article}
\usepackage{hyperref}
\usepackage{amssymb}
\usepackage{amsmath}
\usepackage{amsthm}
\usepackage{mathrsfs}
\newtheorem{theorem}{Theorem}
\newtheorem*{theorema}{Theorem A}
\newtheorem*{theoremb}{Theorem B}

\newtheorem{lemma}{Lemma}
\def\dist{\operatorname{dist}}
\def\re{\operatorname{Re}}
\def\im{\operatorname{Im}}
\def\Ima{\operatorname{Im}}

\def\C{\mathbb{C}}
\def\Z{\mathbb{Z}}
\title{Entire functions with two radially distributed values}
\author{Walter Bergweiler, Alexandre Eremenko\thanks{Supported by NSF grant
DMS-1361836.}$\;$ and Aimo Hinkkanen}
\date{}
\begin{document}
\maketitle
\begin{abstract}
We study entire functions whose zeros and one-points lie on distinct
finite systems of rays. General restrictions on these rays are obtained.
Non-trivial examples of entire functions with zeros and one-points
on different rays are constructed,
using the Stokes phenomenon for second order linear differential equations.

\medskip


\medskip

\end{abstract}

\section{Introduction} \label{sec1}
The zeros of an entire function can be arbitrarily assigned,
but in general one cannot assign the preimages of two values
\cite{Nev}. Since this work of Nevanlinna, various necessary conditions
which the sets of zeros and $1$-points of an entire function must satisfy
were found; see, e.g., \cite{Ozawa,RubelYang,Winkler}. Besides an
intrinsic interest, these conditions are relevant to control theory \cite{Bl,BE,E}.

In this paper we study the simplest setting when the zeros and $1$-points lie on
finitely many rays, or are close to finitely many rays.
The word ``ray'' in this paper
will always mean a ray from the origin. For an entire function~$f$,
we say that a value $a$ is {\em radially distributed} if the set $f^{-1}(a)$
is contained in the union of finitely many rays. 

We begin by recalling some classical results.

\begin{theorema}
{\rm (A. Edrei  \cite{Edr})}
Suppose that all  zeros and
$1$-points of an entire function $f$ are distributed on a finite set of rays, and let
$\omega$ be the smallest angle between these rays. Then the order of $f$
is at most $\pi/\omega$. 
\end{theorema}

The most difficult part of this theorem is the conclusion that the order is
finite. Some special cases of Theorem A under the a priori assumption 
that the order is finite were proved by Bieberbach \cite{Bieber}.

%In the special case that all  zeros and $1$-points are on one ray,
%we can apply this result with $\omega=2\pi$
%and conclude that the order of $f$ is at most $1/2$.
%This result follows -- for functions of finite order -- from a result of 
%Bieberbach \cite[Satz~I$'$, p.~178]{Bieber} proved
%in 1919 that if an entire function $f$ of order $\rho$ omits the values  $0$
%and $1$ in a sector of opening $\omega = \alpha \pi$, then either
%$\rho \leq 1/\alpha=\pi/\omega$   or   $\rho \geq 1/(2-\alpha)$. 
%Bieberbach \cite[Satz~II$''''$, p.~189]{Bieber} also showed that given $\rho$ 
%satisfying $0\leq\rho\leq 1/2$ there 
%exists  an entire function $f$ of order $\rho$ for
%which all zeros and $1$-points are positive.
%For $\rho=1/2$ such a function is given by $f(z)=\cos\sqrt{z}$.

In the 1920s, Biernacki~\cite{Bier1} and Milloux~\cite{Mil}
studied the case that 
all zeros of a transcendental entire function $f$ lie
on some line and all $1$-points
lie on a different line. 
Under certain additional hypotheses, including  in particular 
the assumption that the order is not an integer,
Biernacki showed that the lines must be parallel.
Milloux obtained this result under less restrictive hypotheses,
with additional conclusions on the order of the function.
Both Biernacki and Milloux considered only functions of finite order,
but 
by Edrei's Theorem~A we know now that this is a
consequence of the other hypotheses.

The following corollary from our main results completely
describes the situation 
for the case that the lines intersect.
\begin{theorem} \label{thm1}
Suppose that all zeros of an entire function $f$ lie on
a line $L_1$ and all $1$-points lie on a different
line $L_2$ intersecting $L_1$.
Then $f$ is either of the form $f(z)=e^{az+b}$ or $f(z)=1-e^{az+b}$,
or a polynomial of degree at most~$2$.
\end{theorem}
This result is essentially contained in Milloux's paper.
He did not mention the exceptional cases, perhaps assuming implicitly
that there are indeed both zeros and $1$-points.
In his proof, Milloux omitted a detailed treatment of the case when $L_2$
is orthogonal to $L_1$. 
He also omitted the treatment of the case when $f$ is of the form 
$f(z)=P(z) \prod_{n=1}^{\infty} (1-z/a_n)$ where $P$
is a polynomial and $\sum_n 1/|a_n|<\infty$.  
While these problems can be fixed by considering 
the function $f(z)\overline{{f(-\overline{ z})}}$ in a suitable way,
we include a proof of Theorem~\ref{thm1} below, 
based on different arguments and shorter than that of Milloux.
In fact, we obtain Theorem 1 as a corollary of more general results.

The case of parallel lines was studied
by Baker \cite{Ba} and Kobayashi \cite{Kob} 
who showed that if 
all zeros of a transcendental entire function~$f$ lie on a line $L_1$ and
all $1$-points lie on a different line $L_2$ parallel to~$L_1$,
then $f$ has the form $f(z)=P(e^{az})$ with 
some $a\in\C\backslash \{0\}$ and a polynomial~$P$.

Biernacki returned to the subject in 1929 and
considered the more general situation where the zeros and $1$-points do not lie exactly on certain
lines or rays, but are only close to them.
We say that points $a_n$ (e.g., the zeros of an entire function) accumulate in a direction $\varphi$ if, for every $\varepsilon>0$, the sector $\{z\colon |\arg z - \varphi|<\varepsilon\}$ contains infinitely many points $a_n$. If there is only one such direction $\varphi$ (modulo $2\pi$), we say that the points $a_n$ accumulate in the direction $\varphi$ only.

\begin{theoremb}
{\rm (M. Biernacki \cite[p.~533]{Bier2})}
There does not exist a transcendental entire function $f$ of finite order such that the zeros of $f$ accumulate in a certain direction only and the $1$-points of $f$ accumulate in a different direction only.
\end{theoremb}

We complement the above theorems with several results. First we consider the situation when the zeros lie on a ray and the $1$-points lie on two rays. Note that in view of Theorem~A, these assumptions imply that the function is of finite order.

\begin{theorem} \label{thm0}
Let $f$ be a transcendental entire function whose zeros lie on a ray $L_0$ and whose $1$-points
lie on two rays $L_1$ and $L_{-1}$, each of which is distinct from $L_0$. Suppose that
the numbers of zeros and $1$-points are infinite.
Then $\angle(L_0,L_1)=\angle(L_0,L_{-1})<\pi/2$.
\end{theorem}

It is remarkable that there are non-trivial examples of entire functions
whose zeros lie on the positive ray while all $1$-points lie on two rays
that are not contained in the real line.
As far as we know the functions given in the following theorem 
are the first non-trivial examples of entire functions
with zeros and $1$-points on finitely many rays.

\begin{theorem} \label{thm2}
For every integer $m\geq 3$, there exists an entire function $f$
of order $1/2+1/m$
whose zeros are positive and whose $1$-points lie on the two rays
$\left\{ z \colon \arg z=\pm 2\pi/(m+2)\right\}$.
\end{theorem}

Theorem~\ref{thm1} implies that such functions do not exist for $m=2$, except for the trivial
case where $f$ has no zeros at all and is of the form $f(z)=e^{az}$.

Taking $f(z^n)$ with $f$ as in Theorem~\ref{thm2}
we obtain an entire function whose zeros lie on $n$ rays
and whose $1$-points lie on $2n$ rays distinct from those rays where the zeros lie.

Now we relax the condition that the zeros
and $1$-points are radially distributed.
Let
$$A=\bigcup_{j=1}^n A_j,\quad A_j=\{ te^{i\alpha_j} \colon t\geq 0\},
\quad 0\leq \alpha_1<\ldots<\alpha_n<2\pi,$$
be a finite union of rays.
We say that the $a$-points of an entire function $f$
are {\em close} to the set $A$ if they accumulate only in the directions given by~$A$.
Thus for every  $\varepsilon>0$ 
all but finitely many $a$-points are in the union of sectors 
$$\bigcup_{j=1}^n\{ z \colon |\arg  z-\alpha_j|<\varepsilon\}.$$

Our next result describes the possible configurations of 
finite unions of rays $A$ and $B$ which, apart from the origin,
are pairwise disjoint and have the property that the zeros are
close to $A$ and the $1$-points are close to~$B$.
We will assume that the system of
rays $A\cup B$ is {\em minimal} is the sense that for every ray
$\{ t e^{i\alpha}\colon t\geq0\}$ in $A$ (respectively in $B$)
there is a sequence $(z_k)$ tending to $\infty$ such that
$f(z_k)=0$ (respectively $f(z_k)=1$)
for all $k$, and $\arg z_k\to \alpha$ as $k\to\infty$.

\begin{theorem} \label{thm3}
Let $f$ be a transcendental entire function of order $\rho<\infty$ 
whose zeros are close to $A$ and whose $1$-points are close to~$B$,
with $A\cap B=\{ 0 \} $. Suppose that the system $A\cup B$ is minimal.

Then 
\begin{equation} \label{piomega}
\rho=\frac{\pi}{\omega}>\frac12, 
\end{equation}
where $\omega$ is
the {\em largest} angle between adjacent rays in $A\cup B$,
and there exists a system of rays $C=\bigcup_{j=1}^{2m}C_j \subset A\cup B$,
with $m\geq 1$,
partitioning the plane into $2m$ sectors $S_j$ such that
$\partial S_j=C_j\cup C_{j+1}$ for $1\leq j \leq 2m-1$
and $\partial S_{2m} = C_{2m} \cup C_1$,
with the following properties:

$(i)$ The angle of $S_j$ at $0$ is $\pi/\rho$ when $j$ is even,
and at most $\pi/\rho$ when $j$ is odd. 

$(ii)$ Both boundary rays of an odd sector belong to the same set, $A$ or~$B$.

$(iii)$ There are no rays of $A\cup B$ inside the even sectors.

$(iv)$ If there are rays of $A$ inside an odd sector,
then the boundary rays of this sector belong to~$B$.
If there are rays of $B$ inside an odd sector, then the boundary rays of
this sector belong to~$A$.

$(v)$ If there are no rays of $A\cup B$ in an odd sector,
then its opening angle is $\pi/\rho$.
\end{theorem}
Note that Theorem~B is an immediate consequence of Theorem~\ref{thm3}.

The next result -- whose proof we will only sketch -- shows that 
Theorem~\ref{thm3} is best possible.

\begin{theorem} \label{thm4}
Let $A$ and $B$ be systems of rays, 
satisfying conditions $(i)$--$(v)$ of Theorem~$\ref{thm3}$
for some $\rho\in(0,\infty)$. Then there exists
an entire function of order $\rho$ whose zeros are close to
$A$ and whose $1$-points are close to~$B$.

Moreover, for all finite systems of rays $A$ and $B$ there exists an 
entire function of infinite order whose zeros are close to
$A$ and whose $1$-points are close to~$B$.
\end{theorem}
We note that Theorem~\ref{thm4} shows in particular that the hypothesis that 
$f$ be of finite order is essential in Biernacki's Theorem~B.

We illustrate our results by returning to the case of three  distinct rays, 
which was already discussed in Theorems~\ref{thm0} and~\ref{thm2}.
Now we consider three rays
\[
L_j=\{ te^{ij\alpha}\colon t\geq 0\},\quad j\in\{-1,0,1\},
\]
with $\alpha\in (0,\pi)$. Theorems A, \ref{thm1}, \ref{thm2}  and \ref{thm3} imply the following.

Theorem~\ref{thm2} shows that for certain $\alpha\in(0,\pi/2)$
there exists a transcendental entire function of order $\pi/(2\pi-2\alpha)$
whose zeros lie on $L_0$ while its $1$-points lie on $L_1\cup L_{-1}$.
This result has been extended in \cite{E1} to all $\alpha\in (0,\pi/3]$. 
It remains open whether such functions exist for $\alpha\in(\pi/3,\pi/2)$;
see the discussion at the end of the paper on possible generalizations
of this theorem.

If $\alpha=\pi/2$, then, according to Theorem~\ref{thm1},
there is no transcendental entire function with all zeros on $L_0$
and all $1$-points on $L_1\cup L_{-1}$, unless
it omits $0$ or $1$.
However, the entire function $f(z)=1/\Gamma(-z)$ has zeros on $L_0$
and $1$-points close to the imaginary axis.
This follows from Stirling's formula.

Finally, Theorem~\ref{thm3} implies that if $\alpha\in (\pi/2,\pi)$, then there
is no transcendental entire function of finite order whose zeros are close to $L_0$ and whose $1$-points
are close to $L_1\cup L_{-1}$.

This work was stimulated by questions asked by Gary Gundersen.
We thank him for drawing our attention to 
these problems and for interesting discussions.
We also thank the referee for very valuable comments.

The plan of the paper is the following.
In section~\ref{sec2} we first prove Theorem~\ref{thm3} and then deduce Theorem~\ref{thm1}
and Theorem~\ref{thm0} from it.
The proof of Theorem~\ref{thm4} showing the sharpness of  Theorem~\ref{thm3} is 
then sketched  in section~\ref{proofthm4}.
The proof of Theorem~\ref{thm2} is independent of the rest
and will be given in section~\ref{sec4}.

\section{Proof of Theorems~\ref{thm3}, \ref{thm1} and \ref{thm0}} \label{sec2}
\noindent
\emph{Proof of Theorem}~{\ref{thm3}}.
If the $a$-points of $f$ are close to a finite system of rays,
then evidently $f(z+c)$ has the same property for every $c\in\C$,
with the same rays.
Therefore we may assume without loss of generality that
\begin{equation}\label{B}
f(0)\not\in\{0,1\}.
\end{equation}
We use the standard notation
$$M(r)=M(r,f)=\max_{|z|\leq r}|f(z)|.$$
Let  $(r_k)$ be a sequence tending to $\infty$ with the
property that 
\begin{equation}
\label{A}
\log M(tr_k)=O(\log M(r_k)),\quad k\to\infty,
\end{equation}
for every $t>1$.
Such sequences always exist for functions of finite order.

A sequence $(r_k)$ is called a sequence of {\em P\'olya peaks of
order $\lambda\in [0,\infty)$} for $\log M(r)$, if for every $\varepsilon>0$ we have
\begin{equation}\label{pp}
\log M(tr_k)\leq (1+\varepsilon)t^{\lambda}\log M(r_k),
\quad \varepsilon\leq t\leq 
\varepsilon^{-1},
\end{equation}
when $k$ is large enough.
It is clear that every sequence of P\'olya peaks satisfies~(\ref{A}).
According to a result of Drasin and Shea \cite{DS},
P\'olya peaks of order $\lambda$ exist for all finite
$\lambda\in[\rho_*,\rho^*]$, where
\begin{equation}\label{rho1}
\rho^*=\sup\left\{ p\in {\mathbb R}  \colon  \limsup_{r,t\to\infty}
\frac{\log M(tr)}{t^p\log M(r)}=\infty\right\}
\end{equation}
and
$$\rho_*=\inf\left\{ p \in {\mathbb R}  \colon  \liminf_{r,t\to\infty}
\frac{\log M(tr)}{t^p\log M(r)}=0\right\}.$$
We always have
$$0\leq\rho_*\leq\rho\leq\rho^*\leq\infty,$$
so when $\rho<\infty$, then there exist P\'olya peaks of some (finite) order~$\lambda$.

We refer to \cite[Ch.\ III]{Hor}, \cite[Ch.\ III]{Hor2} and \cite{Ran} 
for the basic results on subharmonic functions used below.
Fixing a sequence $(r_k)$ with the property~(\ref{A}), we consider the two
sequences $(u_k)$ and $(v_k)$ of subharmonic functions in ${\mathbb C}$ given by
\begin{equation}\label{duv}
u_k(z) = \frac{\log|f(r_kz)|}{\log M(r_k)}\quad\mbox{and}\quad
v_k(z) =  \frac{\log|f(r_kz)-1|}{\log M(r_k)}.
\end{equation}
In view of~(\ref{A}), these sequences are bounded from above on every
compact subset of $\C$. 
It follows from~(\ref{B}) that the sequences $u_k(0)$ and $v_k(0)$
tend to~$0$.
According to a well known compactness principle (see, for example,
\cite[Theorems 4.1.8, 4.1.9]{Hor} or 
\cite[Theorems 3.2.12, 3.2.13]{Hor2}),
one can choose a subsequence of $(r_k)$,
which we do without changing notation, such
that the limit  functions  
\begin{equation}\label{lim}
u(z)=\lim_{k\to\infty}\frac{\log|f(r_kz)|}{\log M(r_k)}\quad\mbox{and}\quad
v(z)=\lim_{k\to\infty}\frac{\log|f(r_kz)-1|}{\log M(r_k)}
\end{equation}
exist and are subharmonic in ${\mathbb C}$.  Here the convergence is in the Schwartz
space~$\mathscr{D}'$.
It implies the convergence in $L^1_{\mathrm{loc}}$ and also
the convergence of
the Riesz measures, as the Laplacian is continuous in $\mathscr{D}'$.

The functions $u$ and $v$ are non-zero subharmonic functions in ${\mathbb C}$, and
they have the following properties (we write $u^+ = \max\{u,0\}$):
\vspace{.1in}

$(a)$ $u^+=v^+$.
\vspace{.1in}

$(b)$ $\{ z \colon  u(z)<0\}\cap\{ z \colon  v(z)<0\}=\emptyset$.
\vspace{.1in}

$(c)$ $u$ is harmonic in $\C\backslash A$ and $v$ is harmonic in $\C\backslash B$.
\vspace{.1in}

If $(r_k)$ is a sequence of P\'olya peaks of order $\lambda>0$,
then we have the additional property
\vspace{.1in}

$(d)$ $u(0)=v(0)=0$, and $\max\{u(z),v(z)\}\leq |z|^\lambda$ for all $z\in {\mathbb C}$.
\vspace{.1in}

Properties $(a)$ and $(b)$ are evident. Property $(c)$ holds because the
Laplacian is continuous in $\mathscr{D}^\prime$. Property $(d)$ is a 
consequence of~(\ref{B}) and~(\ref{pp}). Indeed,~(\ref{B}) and 
$$u(0)\geq\limsup_{k\to\infty} u_k(0),$$
(see \cite[(4.1.8)]{Hor}) imply that $u(0)\geq 0$, while~(\ref{pp}) yields
$u(z)\leq|z|^\lambda$ and thus, in particular, $u(0)=0$. The same argument applies
to~$v$.

The components of the complement $\C\backslash(A\cup B)$
will be called {\em sectors of the system} $A\cup B$.

\begin{lemma} \label{le1}
Let $u$ and $v$ be two non-zero subharmonic functions in the plane
which satisfy $(a)$, $(b)$ and $(c)$.  Then either $u(z)\equiv v(z)\equiv c$ for some
$c>0$, or there exist an even number of rays $C_1,\ldots,C_{2m}$, with $m\geq 1$,
that belong to $A\cup B$ and partition the plane into
sectors $S_j$, so that $\partial S_j=C_j\cup C_{j+1}$ for $1\leq j\leq 2m-1$ and
$\partial S_{2m}= C_{2m}\cup C_1$, such that
$u(z)=v(z)>0$ for $z$ in the even sectors
while
$u(z)\leq 0$ and $v(z)\leq 0$ for $z$ in the odd sectors.
Moreover, in each odd sector one of the functions $u,v$ is strictly negative
and the  other is zero.
If $u(z)=0$ in $S_{2k+1}$ then $\partial S_{2k+1}\subset A$,
and if $v(z)=0$ in $S_{2k+1}$ then $\partial S_{2k+1}\subset B$.
\end{lemma}

\noindent
{\em Proof of Lemma}~\ref{le1}. 
If $D$ is a sector of the system $A\cup B$,
and if at some point $z_0$ in~$D$, we have $\max\{ u(z_0),v(z_0)\}>0$,
then $u(z)=v(z)>0$ for
all points $z\in D$. Indeed, both $u$ and $v$ are harmonic in $D$ by $(c)$,
and $(a)$ gives 
\begin{equation}\label{aa}
u(z_0)=v(z_0)>0.
\end{equation}
If $\min\{ u(z_1),v(z_1)\}<0$ for some $z_1\in D$, then this also
holds in a neighborhood of $z_1$, and one of the functions $u$ and $v$ must be zero
in this neighborhood by~$(b)$. Then it is identically equal to zero in $D$
which contradicts~(\ref{aa}). Thus $u$ and $v$ are non-negative in~$D$,
and the minimum principle implies that they are positive. Then they are equal
in $D$ by~$(a)$.
Such sectors $D$ will be called {\em positive} sectors.

If one of the functions $u$ and $v$ is constant, then both functions
are equal to the same positive constant. This follows from $(a)$ and~$(b)$ and the assumption that neither $u$ nor $v$ is identically zero.
For the rest of the proof we assume that both $u$ and $v$ are non-constant.

Suppose that some ray $L\subset A\cup B$ has the property
that positive sectors $D_1$ and $D_2$
are adjacent to $L$
on both sides, that is, $L=\partial D_1\cap\partial D_2$.
(We will see in a moment that $D_1\neq D_2$).
Then we have $u(z)=v(z)$ for $z\in D=D_1\cup D_2\cup L$, in view of~$(a)$,
and $u$ and $v$ must be positive and harmonic in $D$ in view of~$(c)$.

If there are no non-positive sectors, then $u$ and $v$ are equal, positive
and harmonic in $\C\backslash\{0\}$, which is impossible under the current
assumption that they are non-constant. So there
is at least one non-positive sector. In particular, $D_1\neq D_2$ in
the previous paragraph.

Let $D$ be a positive sector. Let $z_0$ be
a point on $\partial D\cap\partial D'$, where $D'$
is a non-positive sector. This means
that $u(z)\leq 0$ and $v(z)\leq 0$ in $D'$.
Then $u(z_0)=v(z_0)=0$. Indeed, $u(z_0) \geq 0$ and $v(z_0)\geq 0$
by the upper semi-continuity of subharmonic functions.
As $D'$ is not
thin at $z_0$ (in the sense of potential theory, see~\cite{Ran})
we obtain that $u(z_0)=v(z_0)=0$. 

Let $C$ be the union of those rays in $A\cup B$ which separate
a positive and a non-positive sector.
It follows from the above considerations that $C$ 
can be written in the
form $C=\bigcup_{j=1}^{2m}C_j$, with $m\geq 1$, with rays $C_j\subset A\cup B$ 
so that in the sector $S_j$ between $C_j$ and $C_{j+1}$ the functions $u$ and $v$
are positive for even $j$ and non-positive for odd~$j$.
Moreover, we have $u(z)=v(z)=0$ for $z\in C$. Note that the sectors with respect to
the system $C$ may be unions of several sectors and rays of the system $A\cup B$. 

Let $S_{2j-1}$ be an odd sector.
Then $u(z)\leq 0$ in $S_{2j-1}$. If $u(z)=0$
in $S_{2j-1}$, then $u$ is not harmonic on either
of the two rays in $\partial S_{2j-1}$,
so $v$ is harmonic on these two rays. As $v(z)>0$ in the adjacent sectors
to $S_{2j-1}$, we obtain that $v(z)\not\equiv 0$ in $S_{2j-1}$, so $v(z)<0$
in $S_{2j-1}$ by maximum principle. In this case we evidently have $\partial S_{2j-1}\subset A$.
Similarly, if $v(z)=0$ in $S_{2j-1}$ then $u(z)<0$ in $S_{2j-1}$ and 
$\partial S_{2j-1}\subset B$.

This proves Lemma~\ref{le1}.
\vspace{.1in}

We return to the proof of Theorem~\ref{thm3}.
Lemma~\ref{le1} does not exclude the possibility that the set of rays $C_j$
is empty, and thus the whole plane coincides with one positive sector.
In this case $u$ and $v$ are identically equal to the same positive constant. 
The following argument shows that this is impossible.

Suppose first that 
\begin{equation}\label{ass}
\rho^*=0
\end{equation}
in (\ref{rho1}). Then
\begin{equation}\label{raz}
\log M(tr)\leq c_0t^{1/4}\log M(r),\quad t>t_0,\; r>r_0,
\end{equation}
with some $c_0>0$.
This implies (\ref{A}), thus from every sequence $r_k\to\infty$ one can choose
a subsequence such that the functions (\ref{duv}) have limits $u,v$ satisfying
Lemma~\ref{le1}.

We claim that these functions must be constant when (\ref{raz}) holds.
Indeed, (\ref{raz}) implies that
\begin{equation}\label{dwa}
u(z)\leq c_0|z|^{1/4},\quad |z|>t_0.
\end{equation}
If $u$ is not constant, then there exists at least one positive sector $D$
of the system $C$, so that $u(z)>0,\; z\in D$ and $u(z)=0,\; z\in\partial D$,
and $u$ is harmonic in $D$. But every positive harmonic function
in a sector, zero on the boundary must have the form
\begin{equation}\label{form}
c|z|^{\pi/\alpha}\cos\left(\frac{\pi}{\alpha}(\arg z-\theta_0)\right),
\quad z\in D,
\end{equation}
where $\alpha\leq 2\pi$ is the opening angle of the sector. This can be seen by 
transforming the sector $D$ conformally to a half-plane, for which the result is standard, see, for example,
\cite[Theorem I]{Boas}.
So $\pi/\alpha\geq 1/2$
in contradiction to (\ref{dwa}). This proves the claim that
$u$ and $v$ coincide with the same positive constant function 
when (\ref{raz}) holds.

Let $z_k=r_ke^{i\beta_k}\to\infty$ be a sequence of zeros of $f$,
such that $\beta_k\to\beta$
where $\beta$ is an argument of some ray of the system $A$.
Using the sequence $r_k$ we define functions $u_k,v_k$ by (\ref{duv}),
and consider some limit functions $u,v$. These functions are both equal
to a positive constant $c$ under our assumption (\ref{ass}).

We have $u_k\to u$ in $\mathscr{D}'$. According to Azarin \cite{A},
this convergence also holds in the following sense: for every $\varepsilon>0$ the set
$$\{ z\colon  |z|<2,\;|u(z)-u_k(z)|>\varepsilon\}$$
can be covered by disks the sum of whose radii is at most $\varepsilon$,
when $k$ is large enough. Let $D$ be the closed disk with the center
at $e^{i\beta}$ of radius $\delta<1/2$, where $\delta$ is also so small that $D$ can meet only one ray in $A\cup B$.
Choosing $\varepsilon<\min\{\delta/2,c/2\}$, we see that for large $k$,
there is a disk $B_k$ around $e^{i\beta}$ such that $u_k(z)\geq u(z)/2=c/2>0$ for $z\in\partial B_k$.
This means that $\log|f(z)|\geq(c/3)\log M(r_k)$ for $z$ on the circles
$$\{ z\colon  z/r_k\in\partial B_k\}.$$
Each of these circles encloses the zero $z_k$ of $f$. Thus by Rouch\'e's theorem,
each of them also contains a $1$-point of $f$. As $\delta$ can be arbitrarily small,
the argument of this $1$-point is close to $\beta$.
This is a contradiction with our
assumption that the zeros and $1$-points are close to disjoint systems of rays.
This contradiction shows that 
\begin{equation}\label{be}
\rho^*>0.
\end{equation}

So for some $\lambda>0$ there exists a sequence of P\'olya peaks
of order $\lambda$. Using this sequence to define $u_k,v_k$ as in (\ref{duv}),
we obtain limit functions $u,v$ satisfying conditions $(a),(b),(c)$ of Lemma~\ref{le1},
and in addition, condition $(d)$. This condition $(d)$
excludes the possibility
that $u$ and $v$ are equal to the same constant in Lemma~\ref{le1}.
So there exist at least one positive and at least one non-positive sector of the system $C$. If $D$ is a positive sector
of opening $\alpha$, then  (\ref{form}) holds in $D$. Comparing this with
the condition $(d)$ before Lemma~\ref{le1}, we conclude that $\pi/\alpha=\lambda$.
As $\alpha$ is an angle between certain two rays of a finite system $A\cup B$,
there are only finitely many possibilities for $\lambda$. On the other hand,
the possible orders $\lambda$ of P\'olya peaks of a given function fill an
interval $[\rho_*,\rho^*]$ which contains $\rho$. We conclude that
$$\rho=\rho_*=\rho^*,$$
and in particular, $\rho^*$ is finite and $\rho>0$.

Using (\ref{rho1}), we obtain 
$$\log M(tr)\leq t^{\rho+1}\log M(r),\quad r>r_0,\; t>t_0.$$
This implies that (\ref{A}) holds for every sequence $r_k\to\infty$.

Also this shows that the angle of every even sector
at the origin is equal to $\pi/\rho$, proving the first statement of $(i)$
of Theorem~\ref{thm3}. 

For $r_0>0$ such that $M(r_0)>1$ 
we consider the curve mapping $[r_0,\infty)$ to $\mathscr{D}'\times\mathscr{D}'$ given by
$$r\mapsto\left(\frac{\log|f(rz)|}{\log M(r)},\frac{\log|f(rz)-1|}{\log M(r)}
\right).$$
Let $F$ be the limit set of this curve when $r\to\infty$. It
consists of pairs $(u,v)$ satisfying $(a)$, $(b)$ and $(c)$ and thus
satisfying the conclusions of  Lemma~\ref{le1}.
As a limit set of a curve, $F$ is closed and connected.
In each sector of the system
$A\cup B$ either both of the functions $u$ and $v$ are positive,
or one is negative.
We conclude that the sectors $S_j$ and their classification into
three classes (positive sectors, those where $u(z)<0$ and those where $v(z)<0$)
is independent on
the choice of the sequence $r_k$.
In particular, $u,v$ are never constant even if the
$r_k$ are not P\'{o}lya peaks.
Further, by Lemma~\ref{le1} we deduce that $(ii)$ of  
Theorem~\ref{thm3} holds.

One consequence of this is that $f(z)\to 0$ in every closed subsector $T$ of an
odd sector where $u(z)<0$. Indeed, if $T$ contains points $z_k$ where $|f(z_k)|\geq \delta>0$ with $z_k\to \infty$ and $\arg z_k\to\beta$,
we choose $r_k=|z_k|$ and, after choosing a subsequence,
consider the limits (\ref{lim}). As $u(e^{i\beta})<0$, we have
$u(z)<-\varepsilon$ in a neighborhood of $e^{i\beta}$ for some $\varepsilon>0$.
 Convergence $u_k\to u$ in
$\mathscr{D'}$ implies convergence in $L^1_{\mathrm{loc}}$,
and we obtain 
$$u_k(w)\leq\frac{1}{\pi r^2}\int_{|z-w|<r}u_k(z)dxdy<-\varepsilon/2,\quad z=x+iy,$$
where $k$ is large enough, $w$ is in a neighborhood of $e^{i\beta}$,
and $r>0$ is small enough. This contradicts our assumption that $|f(z_k)|\geq\delta$ and proves that $f(z)\to 0$ in $T$. 
So such a sector cannot contain rays of $B$.

Similarly, an odd sector in which $v(z)<0$ cannot contain rays of $A$.
This proves statement $(iv)$ of Theorem~\ref{thm3}.

Next we prove statement $(iii)$.
Suppose that a ray $L=\{ te^{i\beta} \colon  t\geq 0\}$
of the set $A$ lies inside an even sector $S_{2j}$.
By assumption, there is  an infinite sequence of zeros $(z_k)$
of the form  $z_k=r_ke^{i\beta_k}$ with $r_k\to\infty$ and $\beta_k\to \beta$.
Passing to a subsequence we may assume that the limits in~(\ref{lim}) exist. 
We obtain a positive harmonic function $u$ in a neighborhood
of $e^{i\beta}$. Using the same argument as before,
with  reference to Azarin's and Rouch\'e's
theorems we conclude that there must be a sequence of $1$-points whose arguments
tend to $\beta$, and this contradicts our assumption that zeros and $1$-points
are close to disjoint systems of rays.

A similar argument shows that there are no
rays from $B$ inside any even sector $S_{2j}$. 

This 
proves $(iii)$ and the fact that the even sectors of the system $C$
coincide with the positive sectors of the system $A\cup B$. 

We next prove 
the second statement of $(i)$. 
Let $S_{2j-1}$ be an odd sector with angle $\pi/\gamma$ at the origin.
Consider again
the limit functions $u$ and $v$ obtained from the P\'olya peaks  $r_k$.
Then we have $(d)$ with $\lambda=\rho$.
One of the two subharmonic functions, say~$u$,
is negative in $S_{2j-1}$, and zero on the boundary of $S_{2j-1}$. Let $h$ be the least
harmonic majorant of $u$
in $G=S_{2j-1}\cap\{ z\colon  |z|<1\}$.
Then $h$ is a negative harmonic function in~$G$, equal to zero on the
straight segments of $\partial G$. It follows that
$$\int_{re^{it}\in G} h(re^{it})\, dt\leq -cr^\gamma, \quad r<1,$$
where $c>0$ and $\gamma=\pi/\alpha$; here $\alpha$ is the angle
of $G$ at the origin. To prove this, notice that the normal derivative of $h$
cannot be zero on the rectilinear part of the boundary of $G$ by the
reflection principle, thus $h$ has a harmonic majorant
of the form
$-cr^\gamma\cos\left(\gamma(t-t_0)\right)$, 
where $t_0$ is the argument of the bisector of $G$.
Then $u$ satisfies the same inequality
$$
u(re^{it})  \leq -cr^\gamma\cos\left(\gamma(t-t_0)\right)  
$$
in $G$. Combining this with property $(d)$ and using, for $0<r<1$, the fact that 
$$0=2\pi u(0)\leq\int_0^{2\pi}  u(re^{it})  \, dt
\leq -c\int_{re^{it} \in G} r^\gamma\cos\left(\gamma(t-t_0)\right) \, dt
+  \int_{re^{it} \notin G}  r^{\rho}\, dt ,
$$
we obtain that $\gamma\geq\rho$ so that $\alpha\leq\pi/\rho$.
This proves the second part of~$(i)$. By this and $(iii)$, we also obtain $\rho=\pi/\omega$ in (\ref{piomega}). That $\rho>1/2$, that is, $\omega<2\pi$, follows since there are at least two sectors.

Finally,
if there are no rays of $A$ inside $S_{2j-1}$, then $u$ is harmonic in
$S_{2j-1}$. Hence, if $u<0$ in $S_{2j-1}$, then $u$ is of the form~(\ref{form}), and it is a harmonic continuation
from an adjacent even sector, so we must have $\alpha=\pi/\rho$. A similar argument
applies if $v$ is negative and there are no rays of $B$ inside $S_{2j-1}$.
This proves $(v)$ and completes the proof of Theorem~\ref{thm3}.

\vspace{.1in}

\noindent
\emph{Proof of Theorem~\ref{thm1}}.
According to Theorem A, the order $\rho$ of $f$ is finite.

First we deal with the case when $f$ is a polynomial,
following Baker \cite{Ba0}.
Without loss of generality, we may assume that $L_1$ is the real line.
Then
$f=cg$, where $g$ is a real polynomial with
all zeros real.
Thus all zeros of $f'$ are real.
Similarly we conclude that all zeros of $f'$ lie on $L_2$,
and hence the point $z_0$ of intersection of $L_1$ and $L_2$
is the only possible zero of $f'$.
Thus $f(z)=c_1(z-z_0)^n+c_2$ for some $n\geq 1$ and
some $c_1,c_2\in {\mathbb C}$ with $c_1\not= 0$.
Such a function $f$ can satisfy the assumptions of Theorem~\ref{thm1} only if 
$f$ is a polynomial of degree at most~$2$.
Notice that this argument can be extended to functions of order less than~$2$,
but we do not use this.

Suppose now that $f$ is transcendental. 
Then we use Theorem~\ref{thm3}.
This theorem implies that there exists
at least one even sector. If there is only one even sector,
and its angle is greater than $\pi$,
then the odd sector does not contain
rays of $A\cup B$, so by $(v)$ its opening must be the same
as the opening of the even sector, which is a contradiction.

If there are two even sectors, then the odd sectors contain no rays of $A\cup B$.
It follows from $(i)$ and $(v)$ that all sectors must have opening $\pi/2$. Then by $(ii)$ the
zeros are close to the boundary of a quadrant, and the $1$-points are close to
the the boundary of the opposite quadrant. But by hypothesis 
the zeros lie on a line and
the $1$-points lie on a line. We conclude that the zeros are actually close to
one ray and the $1$-points are close to another ray.
But then there is only one even sector, contradicting
our assumption at the beginning of this paragraph.

The only remaining possibility is that
there is one even sector with opening $\pi$.
Then $\rho=1$, and we assume without loss of generality
that this sector is the upper half plane and the $1$-points are real.
This means that $B$ consists of two rays whose union is the real line.
Using the notation of the proof of Theorem~\ref{thm3}, and choosing a 
sequence  $(r_k)$ of P\'olya peaks of order~$1$, 
we obtain $u(z) =  \Ima z$ and $v(z) =  \Ima^+ z$.
This implies that 
\begin{equation}\label{N}
N(r_k,1,f)\sim \frac{1}{\pi}\log M(r_k),\quad k\to\infty.
\end{equation}
Now let $g(z)=f(z)\overline{f(\overline{z})}$.
As all $1$-points of $f$ are real, $f(z)=1$ implies that
$g(z)=1$, so if $g\not\equiv 1$, we will have
from~(\ref{N}) that
\begin{equation}\label{Ng}
N(r_k,1,g)\geq (1-o(1) ) \frac{1}{\pi} \log M(r_k),\quad k\to\infty.
\end{equation}
Now define the subharmonic function $w$ in the plane by 
$$w(z)=\lim_{k\to\infty}\frac{\log|g(r_kz)|}{\log M(r_k)}.$$
Above, $M(r_k)$ always refers to $M(r_k,f)$. 
It is evident that $w(z)=u(z)+u(\overline{z})=0$.
Together with~(\ref{Ng}) this implies that
$g(z)\equiv 1$. We conclude that with this normalization,
$f$ has the form $f(z)=\exp(icz+id)$, where  $c$ and $d$ are real.
This completes the proof of Theorem~\ref{thm1}.

\vspace{.1in}

\noindent
\emph{Proof of Theorem~\ref{thm0}}.
We assume without loss of generality that $L_0$ is the
positive ray. The order of $f$ must be finite by Theorem A, so Theorem~\ref{thm3} is
applicable. As there are only three rays, the number $m$ in Theorem~\ref{thm3} must
be~$1$. So we have one even sector of opening $\pi/\rho$ and one odd
sector of opening at most $\pi/\rho$. Thus $\rho\leq 1$. 
In view of~$(ii)$,
the common boundary of the odd and even sector is $L_1\cup L_{-1}$. 
So $L_0$ lies inside the odd sector.
The possibility
that $\rho=1$ is excluded by~(\ref{piomega}) and  Theorem~1.
Hence $\rho<1$,  and the function $f$ is of genus zero and thus of the form
$$f(z)=cg(z),\quad g(z)=z^n\prod_{j=1}^\infty\left(1-\frac{z}{z_k}\right),$$
where $c\in {\mathbb C}\setminus \{0\}$ and $n$ is a non-negative integer, and $(z_k)$ is a sequence of positive numbers
tending to~$\infty$. If $c$ is real,
we conclude that the rays $L_1$ and $L_{-1}$ are symmetric with respect
to $L_0$ which proves Theorem~\ref{thm0} in this case.

Suppose now that $L_{1}$ and $L_{-1}$ are not symmetric with
respect to $L_0$ so that $c$ is not real.
Let us set $a=1/c$. Then $f(z)=1$ is equivalent to $g(z)=a$. We consider 
the function $h(z)=(g(z)-a)/(\overline{a}-a)$. In view of the symmetry of~$g$,
the zeros of $h$ lie on the rays $L_1$ and $L_{-1}$, while the $1$-points lie on
the reflected rays $\overline{L_1}$ and $\overline{L_{-1}}$.
Since $L_0$ lies in the odd sector, which has angle $<\pi$ at the origin, it follows that the two rays $L_1$, $L_{-1}$ are interlaced
with the two rays $\overline{L_1}$, $\overline{L_{-1}}$. 
This contradicts $(ii)$ and $(iii)$ of Theorem~\ref{thm3} applied to $h$, and completes the proof
of Theorem~\ref{thm0}.

\section{Sketch of the proof of Theorem~\ref{thm4}} \label{proofthm4}
We only indicate the construction of examples showing that Theorem~\ref{thm3}
is best possible, as this construction is well-known,
see for example \cite{Drasin},
where a similar construction was used for the first time.

We fix $\rho\in (1/2,\infty)$
and construct a $\rho$-trigonometrically convex function $h$ such that the
union of the even sectors coincides with the set 
$$\{ re^{it}\colon r> 0,\; h(t)>0\},$$
and such that $h$ is trigonometric except at the arguments of some rays
inside the odd sectors. If there are no rays of $A\cup B$
in the odd sectors at all, then $\rho$ must be an integer, and we 
just take $h$ to be of the form $h(t)=\cos ( \rho ( t - t_0 )  )$, where
$t_0$ is the argument of the bisector of one of the sectors.

Then we discretize the Riesz mass of the subharmonic function
$$w(re^{it})=r^\rho h(t),$$
as it is done in \cite{A},
and obtain an entire function $g$ with zeros on some rays $A\cup B$ which
lie in the odd sectors, and such that 
$$\lim_{r \to\infty}r^{-\rho}\log|g(rz)|=w(z).$$
If there are odd sectors with opening $\pi/\rho$, then $h$ must be
trigonometric on the intervals corresponding to these sectors,
so we multiply $g$ by a canonical product of order smaller than $\rho$
to achieve that $g$ has infinitely many zeros on all those
rays of the system $A\cup B$ which lie inside the odd sectors.
Then we label the odd sectors with labels $0$ and $1$:
if the boundary of an odd sector belongs to~$A$, we label it with~$1$,
and if the boundary belongs to $B$ we label it with~$0$.

Let $S_j$ be an odd sector labeled with~$1$. Consider the component $D_j$
of the set 
$\{ z \colon  |g(z)|<2\}$ which is asymptotic to $S_j$.
Let $p$ be a quasiconformal
map of the disk $\{ z \colon  |z|<2\}$ onto itself,
equal to the identity mapping on the boundary,
and such that $p(0)=1$,  whose complex dilatation is supported on
the set $\{ z \colon  3/2\leq |z|\leq 2\}$. We define
$$G(z)=\left\{\begin{array}{ll} p(g(z)),& z\in\bigcup_j D_j,\\
g(z),&\mbox{otherwise}.\end{array}\right.$$
Here the union is over all odd sectors labeled with~$1$.
This $G$ is a quasiregular map of the plane, whose complex dilatation is supported
on a small set $E$ in the sense that
$$\int_E\frac{dxdy}{x^2+y^2}<\infty.$$
Then the theorem of Teichm\"uller--Wittich--Belinski~\cite[\S V.6]{LV}
guarantees the existence
of a quasiconformal map $\phi$ such that $f=G\circ\phi$ is an entire function,
and $\phi(z)/z \to 1$ as $z\to\infty$. It is easy to verify that $f$ has
all the required properties. 
\vspace{.1in}

For the construction of infinite order functions,
let $A=\bigcup_{j=1}^m \{te^{i\alpha_j}\colon t\geq 0\}$
and $B=\bigcup_{k=1}^n \{te^{i\beta_k}\colon t\geq 0\}$ be two 
finite  systems of rays with $A\cap B=\{0\}$.
Again we only sketch the argument.

First we note that by~\cite[Part III, Problems 158--160]{PolyaSzego} there exists 
an entire function $E$ such that $z^2(E(z)+1/z)$ is bounded outside the
half-strip $S=\{z\colon \re z> 0,|\im z|< \pi\}$.
In particular, $E$ is bounded outside~$S$.
Considering $F(z)=\delta(E(z)-c)/((z-a)(z-b))$, where $\delta>0$ is small,
 $c\in\C$, and $a$ and $b$ are $c$-points of~$E$, we obtain an entire 
function $F$ such that 
$$
|F(z)|\leq \frac{1}{|z|^2 + 1}
\leq \frac{1}{\dist(z,S)^2}, \quad z\notin S,
$$
where $\dist(z,S)$ denotes the distance from $z$ to~$S$.
For some large $R>0$ we now consider the functions
$$
a_j(z)=1+z\exp F(e^{-i\alpha_j}z-R)
\quad\text{and}\quad
b_k(z)=z\exp F(e^{-i\beta_k}z-R).
$$
With $S_j=\{e^{i\alpha_j}(z+R)\colon z\in S\}$
we find, noting that $|e^w-1|\leq 2|w|$ for $|w|\leq 1$, that
\begin{align*}
|a_j(z)-z-1|&\leq \left|z\left(\exp F(e^{-i\alpha_j}z-R)-1\right)\right|
\leq 2|z F(e^{-i\alpha_j}z-R)|
\\ &
\leq  \frac{2|z|}{\dist(e^{-i\alpha_j}z-R,S)^2}
=  \frac{2|z|}{\dist(z,S_j)^2}, \quad z\notin S_j.
\end{align*}
Similarly, with 
$T_k=\{e^{i\beta_k}(z+R)\colon z\in S\}$ we have
$$
|b_k(z)-z|
\leq  \frac{2|z|}{\dist(z,T_k)^2}, \quad z\notin T_k.
$$
We choose $\varepsilon>0$ so small that the sectors
$U_j=\{z\colon |\arg(z- e^{i \alpha_j} R/2)-\alpha_j|\leq \varepsilon\}$ and
$V_k=\{z\colon |\arg(z-  e^{i \beta_k}  R/2)-\beta_k|\leq \varepsilon\}$ are disjoint and put
$$
G(z)=
\begin{cases}
a_j(z), \quad z\in U_j,\\
b_k(z), \quad z\in V_k.
\end{cases}
$$
Then 
$$G(z)=z+O(1), \quad z\in \bigcup_{j=1}^m\partial U_j \cup \bigcup_{k=1}^n\partial V_k .$$
This allows us to extend $G$ to a quasiregular map of the plane 
which satisfies 
$$G(z)=z+O(1), \quad z\in \C\backslash \left(\bigcup_{j=1}^m U_j \cup \bigcup_{k=1}^n V_k \right)$$
and whose dilatation $K_G$ satisfies
$K_G(z)=1+O(1/|z|)$ as $z\to\infty$.
Again  the theorem of Teichm\"uller--Wittich--Belinski
yields the existence of a quasiconformal map $\phi$ such that $f=G\circ\phi$ is entire
and $\phi(z)/z \to 1$  as $z\to\infty$. It is not difficult to show that   
the zeros of $f$ are close to $A$ and the $1$-points of $f$ are close to~$B$.

We note that the method does not actually require that the rays that form
$A$ are distinct from those that form~$B$. Indeed, if we want that both zeros and
$1$-points accumulate at $\{te^{i\alpha_j}\colon t\geq 0\}$, we only have to choose
$a_j(z)=c+z\exp F(e^{-i\alpha_j}z-R)$ with a constant $c$  different from $0$ and~$1$.

\section{Proof of Theorem~\ref{thm2}} \label{sec4}
We consider differential equations
\begin{equation}
\label{1}
-y^{\prime\prime}+ \left((-1)^\ell z^m+E\right)y=0, \quad\ell\in\{0,1\},
\quad m\geq 3,
\quad E\in {\mathbb C}.
\end{equation}
Here $m$ is an integer, so all solutions are entire functions.
The equation has the following symmetry property. Set
$$\varepsilon=e^{\pi i/(m+2)},\quad \omega=\varepsilon^2.$$
If $y_0(z,E)$ is a solution of~(\ref{1}) then
\begin{equation}\label{df}
y_k(z,E)=y_0(\omega^{-k}z,\omega^{2k}E)
\end{equation}
satisfies the same equation, while
$$y_0(\varepsilon^{-k}z,\varepsilon^{2k}E)$$
with an odd $k$ satisfies~(\ref{1}) with $(-1)^\ell z^m$ replaced by $(-1)^{\ell +1} z^m$.

The Stokes sectors are defined as follows.
When $\ell=0$, they are $S_0=\{ z \colon  |\arg z|<\pi/(m+2)\}$,
and $S_k=\omega^kS_0$ for $k\in\Z$.  When $\ell=1$, the Stokes sectors are
$S_0=\{ z \colon  0< \arg z <2\pi/(m+2)\}$ and $S_k=\omega^kS_0$ for $k\in\Z$.

To obtain a discrete sequence of eigenvalues, one imposes  boundary
conditions of the form
\begin{equation}
\label{bc}
y(z)\to 0,\quad z\to\infty,\quad z\in S_n\cup S_k,
\end{equation}
for some $n$ and~$k$. The exact meaning of~(\ref{bc}) is that $y(z)\to 0$
when $z\to\infty$ along any interior ray from the origin contained
in the union of the two sectors. 

We will denote such a boundary condition by $(n,k)$. It is known~\cite{S} that
when $n\neq k\pm1$ (modulo $m+2$), then the boundary value problem $(n,k)$ has a discrete
spectrum with a sequence of eigenvalues tending to infinity.
(For completeness, we include the argument below.)
Moreover, K. Shin \cite{Shin2} proved that these eigenvalues always
lie on a ray from the origin. In particular, when $S_n$ and $S_k$
are symmetric with respect to the positive ray, these eigenvalues
are positive. All other cases can be reduced to this case using the symmetry
of the differential equation stated above: if $\omega_1$ and $\omega_2$
are bisectors of $S_n$ and $S_k$, then the eigenvalues lie on the
ray $\{ t/(\omega_1\omega_2)\colon  t\geq 0\}$.

From now on we assume that $\ell=0$ in~(\ref{1}).
For each $E$ the equation~(\ref{1}) has a solution tending to zero as $z\to\infty$ in $S_0$.
More precisely, there is a unique solution $y_0(z,E)$ satisfying
\begin{equation} \label{as}
y_0(z,E)=(1+o(1))z^{-m/4}\exp\left(-\frac{2}{m+2}z^{(m+2)/2}\right)
\end{equation}
as $z\to\infty$ in any closed subsector of $S_0\cup S_1\cup S_{-1}$;
see \cite[Thm 6.1]{S}.
Notice the simple but important fact that this principal part of the
asymptotics does not depend on~$E$.
The function $y_0(z,E)$  is actually an entire function of the two variables $z$ and $E$,
and its asymptotics when $E\to\infty$ while $z$ is fixed are also known
\cite[Thm 19.1]{S};
this implies that the entire function $E\mapsto y_0(z_0,E)$ has order
$$\rho=\frac{1}{2}+\frac{1}{m}.$$

Now we define $y_k$ by~(\ref{df}).
Then $y_k\to 0$ as $z\to\infty$ in $S_k$.
The boundary problem $(n,k)$ thus has a solution when
$y_n$ and $y_k$ are linearly dependent as functions of~$z$.
This means that their Wronskian vanishes.
But the Wronskian, evaluated at $z=0$, is an entire function of $E$,
and its order is less than~$1$. Thus
its zeros, which are the eigenvalues of the problem, form
a sequence tending to infinity, as mentioned above.

As $y_0,y_1,y_{-1}$ satisfy the same differential equation, we have
$$y_{-1}=C(E)y_0+\tilde{C}(E)y_{1}.$$
The asymptotics of $y_{1}$ and $y_{-1}$ in $S_0$ (which follow from~(\ref{as}))
show that 
$\tilde{C}=-\omega$, so
\begin{equation}\label{2}
y_{-1}=C(E)y_0-\omega y_{1}.
\end{equation}
By differentiating this with respect to $z$ we obtain
\begin{equation}
\label{2'}
y_{-1}^\prime=C(E)y_0^\prime-\omega y_{1}^\prime.
\end{equation}
Solving~(\ref{2}) and~(\ref{2'}) by Cramer's rule, we obtain
$$C(E)=W_{-1,1}/W_{0,1},$$
where $W_{i,j}$ is the Wronskian of $y_i$ and $y_j$.
This shows that $C$ is an entire function (because $W_{0,1}$ is never $0$).
It has the same order $\rho$ that $y_0$ has as a function of~$E$.

In view of~(\ref{2}), the zeros of $C$ are
exactly the eigenvalues $\lambda_j$ of the
problem~(\ref{bc}) with $(n,k)=(-1,1)$. So all zeros of $C$ are
positive by Shin's result. 
Substituting $(z,E)\mapsto (\omega^{-1}z,\omega^2E)$ to~(\ref{2}), we obtain
$$y_0=C(\omega^2E)y_1-\omega y_2.$$
Using this to eliminate $y_0$ from~(\ref{2}) we obtain
$$y_{-1}=\left(C(E)C(\omega^2E)-\omega\right)y_1-C(E)\omega y_2.$$
We conclude that the zeros of the entire function
\begin{equation}\label{sib}
g(E):=C(E)C(\omega^2E)-\omega
\end{equation}
are the eigenvalues of the problem $(-1,2)$.
Therefore, these zeros lie on the ray
$\{ z=t\omega^{-1}\colon t\geq 0\}$. So if we define
$f(E)=-\omega^{-1} g(\omega^{-1}E)$ and $h(E)=C(E)/\sqrt{\omega}$,
then
$$f(E)=1-h(\omega^{-1}E)h(\omega E),$$
the zeros of $f$ are on the positive ray and the $1$-points on two other rays. 
This completes the proof of Theorem~\ref{thm2}.
\vspace{.1in}

\noindent
{\em Remarks.} Once it is known that two entire functions $C$ and  $g$
satisfy~(\ref{sib}) and the zeros of each function lie on a ray,
the order of both functions and the angles between the rays can be determined
from Theorem~\ref{thm3} and Theorem~\ref{thm0}. 

Equations of the type~(\ref{sib})
occur for the first time in the work of Sibuya
and his students \cite{S,S2,S3} for the simplest case when $m=3$.

It was later discovered that these equations also arise in the 
context of exactly solvable models of statistical mechanics on
two-dimensional lattices and in quantum field theory \cite{DDT1,DDT2}.

The interesting question is to which angles Theorem~\ref{thm2} generalizes.
Concerning the approach that makes use of differential equations, the following comments are in order. 
If $m>2$ is not an integer, equation (\ref{1}) and its solutions are defined
on the Riemann surface of the logarithm, but Sibuya's solution $y_0$
is still entire as a function of $E$.
We found no source where this fact is proved,
but it is stated and used in \cite[p. 576]{DDT1}, \cite[p. R231]{DDT2}
and \cite{Tabara}. Shin's result, which we used above seems to generalize
to non-integer $m\geq 4$, see \cite[Theorem 11]{Shin2} which
we use with $\ell=1$ and $\ell=2$. On the other hand, numerical evidence
in \cite{Bender} (see Figs. 14, 15, 20)
shows that for $m<4$ our function $g(E)$ in (\ref{sib})
does not have radially distributed zeros on one ray,
even if finitely many of the zeros are discarded. 

By extending our method Theorem 3 is generalized in \cite{E1}
to all angles in $(0,\pi/3]$.
The question remains open for angles in $(\pi/3,\pi/2)$
other than $2\pi/5$.

\section*{Correction}
After this paper was published [Math. Proc. Cambridge Philos. Soc.,
\url{https://doi.org/10.1017/S0305004117000305}]
we noticed a small error.
We claimed towards the end of the introduction
that the $1$-points of  the function $f(z)=1/\Gamma(-z)$ 
are close to the imaginary axis. However, this function also has 
infinitely many $1$-points on the positive real axis. 

To construct an example $f$ with positive zeros and $1$-points close to 
the imaginary axis 
we put $x_k=k(\log k)^2$ for $k\geq 2$ and consider the infinite product
\[
f(z)=\prod_{k=2}^\infty \left(1-\frac{z}{x_k}\right).
\]
It follows from~\cite[Chapter 2, \S 5, $(5.26_2)$, p.~81]{GO} that
given $\delta\in (0,\pi/2)$ we have
%there exists $c>0$ such that
\[
\log|f(re^{i\theta})|\sim -r (\log r)^{-1}\cos\theta
\]
for $\delta<\theta<2\pi-\delta$ as $r\to\infty$.
Together with the Phragm\'en-Lindel\"of principle this yields that 
$f(re^{i\theta})\to 0$ as $r\to\infty$ if $|\theta|\leq\delta$.
%if $re^{i\theta}$ avoids small disks around the $x_k$.
Denoting by $z_k=r_ke^{i\theta_k}$ the $1$-points of $f$ we conclude that 
$\cos\theta_k\to 0$ as $k\to\infty$. Thus the $1$-points of $f$ accumulate only
at the imaginary axis.

{\em W. B.: Mathematisches Seminar

Christian-Albrechts-Universit\"at zu Kiel

Ludewig-Meyn-Str. 4

24098 Kiel

Germany

%\newpage
\vspace{.1in}

A. E.: Department of Mathematics

Purdue University

West Lafayette, IN 47907

USA
\vspace{.1in}

A. H.: Department of Mathematics

University of Illinois at Urbana--Champaign

1409 W. Green St.

Urbana, IL 61801

USA}

\end{document}